\documentclass[10pt, reqno]{amsart}
\usepackage{amssymb, amsthm, amsmath, amsfonts}
\usepackage{graphics}
\usepackage[hypertexnames=false]{hyperref}
\usepackage{bbm}
\usepackage[all]{xy}
\usepackage{enumerate}
\usepackage[mathscr]{eucal}
\usepackage{thmtools}
\usepackage{xcolor}

\usepackage{extarrows}
\usepackage{ytableau}
\usepackage{enumitem}
\usepackage{cite}
\usepackage[english]{babel}
\theoremstyle{plain}
\newtheorem{theorem}{Theorem}[section]
\newtheorem{lemma}[theorem]{Lemma}
\newtheorem{proposition}[theorem]{Proposition}

\newtheorem{corollary}[theorem]{Corollary}
\numberwithin{equation}{section}

\theoremstyle{definition}

\newtheorem{definition}[theorem]{Definition}
\newtheorem{example}[theorem]{Example}
\newtheorem{remark}[theorem]{Remark}

\newcommand{\End}{\operatorname{End}}
\newcommand{\Hom}{\operatorname{Hom}}
\newcommand{\Ext}{\operatorname{Ext}}
\newcommand{\Aut}{\operatorname{Aut}}
\newcommand{\Inj}{\operatorname{Inj}}
\newcommand{\Sur}{\operatorname{Sur}}
\newcommand{\Soc}{\operatorname{Soc}}
\newcommand{\Top}{\operatorname{Top}}
\newcommand{\Homeo}{\operatorname{Homeo}}
\newcommand{\Sh}{\operatorname{Sh}}
\newcommand{\Mod}{{\textrm{-}\mathrm{Mod}}}
\newcommand{\Cont}{\operatorname{Cont}}
\newcommand{\UC}{\operatorname{UC}}
\newcommand{\Fin}{\operatorname{Fin}}

\newcommand{\id}{\operatorname{id}}
\newcommand{\op}{\mathrm{op}}

\newcommand{\FA}{\mathrm{FA}}

\newcommand{\M}{\mathscr M}
\newcommand{\E}{\mathscr E}

\newcommand{\C}{\mathscr C}
\newcommand{\D}{\mathscr D}
\newcommand{\K}{\mathcal K}
\newcommand{\T}{\mathcal T}
\newcommand{\Rcal}{\mathcal R}

\title[Continuous Representations of Topological Monoids]{Continuous representations of topological monoids via reconstruction}

\author[Z.X. Di]{Zhenxing Di}
\address{Z.X. Di \ School of Mathematical Sciences, Huaqiao University, Quanzhou 362021, China}
\email{dizhenxing@163.com}
	
\author[L.P. Li]{Liping Li}
\address{L.P. Li \ Department of Mathematics, Hunan Normal University, Changsha 410081, China}
\email{lipingli@hunnu.edu.cn}
	
\author[L. Liang]{Li Liang}
\address{Li Liang: Department of Mathematics, Gansu Provincial Basic Discipline Research Center for Complex Systems Analysis and Control, Lanzhou Jiaotong University, Lanzhou 730070, China}
\email{lliangnju@gmail.com}
\urladdr{https://sites.google.com/site/lliangnju}
	
\thanks{Z. Di was partly supported by NSF of China (Grant No. 12471034), Scientific Research Fund of Fujian Province (Grant No. 605-52525002) and Scientific Research Fund of Huaqiao University (Grant No. 605-50Y22050); L. Li was partly supported by NSF of China (Grant No. 12571037); L. Liang was partly supported by NSF of China (Grant No. 12271230) and NSF of Gansu Province (Grant No. 26JRRA059).}

\begin{document}

\begin{abstract}
We study continuous representations of topological monoids on discrete modules via a finite-restriction reconstruction method. For an arbitrary topological monoid $\M$, we construct a transformation monoid $\E_\M$ whose uniformly continuous representations are equivalent to continuous representations of $\M$ on discrete modules. By the finite-restriction reconstruction in \cite{LiRel}, these categories are equivalent to categories of sheaves of modules over special sites. We illustrate this method by several examples, including various realizations of $\FA$-modules and continuous representations of $\operatorname{End}_A(\varinjlim_r A^r)$ for a finite ring $A$.
\end{abstract}

\maketitle

\section{Introduction}

\subsection{Motivation}

Artin's reconstruction theorem identifies continuous actions of topological groups on discrete sets and sheaf categories on orbit categories endowed with the atomic topology \cite{MacLaneMoerdijk, Jo}. Recently, the same idea has been used to explore continuous actions of topological monoids. For instance, Hemelaer and Rogers \cite{HR} studied relations between algebraic properties of a discrete monoid and topos-theoretic properties of its category of actions; Rogers \cite{Rogers} proved that continuous actions of an arbitrary topological monoid on discrete sets form a Grothendieck topos and developed a Morita theory for these action topoi; he \cite{Rogers1} also established a Galois correspondence between endomorphism monoids and relational extensions of models.

These results motivated us to apply methods from the sheaf theory to investigate continuous representations of topological monoids on modules with discrete topology, extending our work for topological groups in \cite{DLLX}. In \cite{LiRel} we obtained a first step for transformation monoids and uniformly continuous representations. Explicitly, let $\Omega$ be a discrete set and let $\E \leqslant \Omega^{\Omega}$ be a transformation monoid equipped with the pointwise-convergence topology. We introduced the notion of uniformly continuous $k\E$-module and established an of Artin's reconstruction theorem for monoids:
\[
k\E \Mod^{\mathrm{uc}} \simeq \Sh(\C^{\op}, \, J, \, \underline{k}),
\]
where $\C$ is a finite-restriction category with finite subsets of $\Omega$ as objects and restrictions of elements of $\E$ as morphisms, and $J$ is the Grothendieck topology generated by inclusions of finite subsets. For details, see \cite[Subsection 7.4]{LiRel}.

The main purpose of this paper is to extend this sheaf-theoretic reconstruction to continuous representations of arbitrary topological monoids on discrete modules, removing both the transformation-monoid assumption and
the uniform-continuity requirement.

\subsection{Main work}

Given a topological monoid $\M$, we construct a discrete $\M$-set $\Omega_{\M}$ and a transformation monoid $\E_{\M} \leqslant (\Omega_{\M})^{\Omega_\M}$. Our first main result (Theorem \ref{thm:universal-realization}) asserts that
\[
k\M\Mod^{\mathrm{cont}} \simeq k\E_{\M} \Mod^{\mathrm{uc}} = k\E_{\M} \Mod^{\mathrm{cont}}.
\]
Combining these equivalences with the monoidal version of Artin's theorem of \cite{LiRel} and the torsion-theoretic description of sheaves in \cite{DLL}, we obtain the following reconstruction method
\[
k\M\Mod^{\mathrm{cont}} \simeq k\E_{\M}\Mod^{\mathrm{uc}} \simeq \Sh(\C_{\M}^{\op}, \, J, \, \underline{k}) \simeq \C_{\M}\Mod/\T_J,
\]
where $\C_{\M}$ is the finite-restriction category associated with $\E_{\M}$ and $\T_J$ denotes the Serre localizing subcategory of $J$--torsion $\C_{\M}$-modules.

For many concrete examples, the above method can be significantly simplified, as explained below:
\begin{enumerate}
\item The universal construction $\Omega_{\M}$ is often unnecessary. For instance, suppose that $\M \leqslant \Omega^{\Omega}$ is a transformation monoid, and for a cofinal family of finite subsets $\Gamma \subseteq \Omega$, whenever $f, g \in \M$ agree on $\Gamma$ there exist elements $h_1, h_2$ fixing $\Gamma$ pointwise such that $f \circ h_1 = g \circ h_2$. Then every continuous $\M$-action is uniformly continuous. In this situation the first step of the general construction can be omitted completely.

\item The second simplification occurs at the level of the Grothendieck topology. If every finite inclusion in $\C_\M$ admits a retraction, then the Grothendieck topology $J$ is trivial, and the sheaf category is simply a presheaf category. More generally, when $J$ is a rigid subcategory topology, the sheaf category is equivalent to the representation category of the corresponding rigid subcategory.

\item Finally, for several important examples the category $\C_\M$ carries a generalized Reedy structure. One may then apply the methods developed in \cite{LiHom, LiNorm} to construct Dold--Kan-type equivalences, replacing $\C_\M$ by a substantially simpler $k$-linear category.
\end{enumerate}

These reductions can be illustrated by the following familiar transformation monoids.
\begin{itemize}
\item Let $\Omega$ be an infinite discrete set and $\M = \Omega^{\Omega}$. Via the above reductions, one has
    \[
    k\M \Mod^{\mathrm{cont}} = k\M \Mod^{\mathrm{uc}} \simeq \FA \Mod,
    \]
where $\FA$ is the category of nonempty finite sets. Furthermore, by \cite{LiRel} we obtain a classification of irreducible objects of these representation categories over any field.

\item Let $\M = \End(\mathbb{Q}, \leqslant)$. Similarly, one has:
    \[
    k\M\Mod^{\mathrm{cont}} = k\M\Mod^{\mathrm{uc}} \simeq \mathrm{OA} \Mod,
    \]
where $\mathrm{OA}$ is equivalent to the simplex category. By Dold-Kan correspondence, these representation categories are equivalent to the category of nonnegatively graded cochain complexes of $k$-modules.

\item Let $\mathbb{F}_q$ be a finite field, $\Omega = \varinjlim_r \mathbb{F}_q^r$, and $\M = \End_{\mathbb{F}_q} (\Omega)$. We have:
    \[
    k\M\Mod^{\mathrm{cont}} = k\M\Mod^{\mathrm{uc}} \simeq \mathrm{VA} \Mod,
    \]
where $\mathrm{VA}$ is the category of finite dimensional vector spaces and linear maps over $\mathbb{F}_q$. Moreover, when $q$ is invertible in $k$, Kuhn's result \cite{Kuhn} asserts that $k\M \Mod^{\mathrm{cont}}$ is equivalent to the representation category of the groupoid of finite dimensional vector spaces and invertible linear maps.
\end{itemize}

The reconstruction method also explains the following invariance phenomena.
\begin{enumerate}
\item Let $\mathcal{A}$ be an inclusion-cofinal system of objects in $\C_\M$, and let $\D$ be the full subcategory of $\C_\M$ consisting of these objects. Then the comparison lemma tells us that
\[
\Sh(\C_\M^{\op}, \, J, \, \underline{k}) \simeq \Sh(\D^{\op}, \, J_{\D}, \, \underline{k}),
\]
where $J_{\D}$ is the induced Grothendieck topology on $\D^{\op}$. In particular, when $J$ is the trivial topology, one has $\C_\M \Mod \simeq \D \Mod$. Therefore, any full subcategory of $\mathrm{FA}$ (resp., $\mathrm{OA}$, $\mathrm{VA}$, and many other examples) containing infinitely many objects (up to isomorphism) has the same representation theory as the original category.

\item Different endomorphism monoids can have the same continuous representation theory. Let $\Omega$ be an infinite discrete set, $\mathscr{F}=\Omega^\Omega$ the full transformation monoid equipped with the topology of pointwise convergence, and $\M$ a dense submonoid of $\mathscr{F}$. We show that
\[
k\M \Mod^{\mathrm{uc}} \simeq k\mathscr{F} \Mod^{\mathrm{uc}} \simeq \FA \Mod.
\]
In particular, if every continuous $k\M$-module is uniformly continuous, then
\[
k\M \Mod^{\mathrm{cont}} \simeq \FA \Mod.
\]
This applies to several familiar transformation monoids, including the monoid of all surjective self-maps of $\Omega$, as well as the monoids of continuous self-maps of the real line and the Cantor space, equipped with the subspace topology inherited from $\mathscr{F}$ rather than the usual compact-open topology.
\end{enumerate}

As another example, we consider the endomorphism monoids of infinite free modules over finite rings. Explicitly, let $A$ be a finite ring, $\Omega = A^\infty = \varinjlim_r A^r$ with the discrete topology, and $\M = \End_A(\Omega)$ carrying the pointwise-convergence topology. Using the standard coordinate submodules $A^r \subseteq \Omega$ whose inclusions split, we obtain
\[
k\M \Mod^{\mathrm{cont}} \simeq \D_A \Mod,
\]
where $\D_A$ is the category of finite free $A$-modules and all $A$-linear maps. The above equivalence realizes $\D_A \Mod$ as the category of continuous representations of a single infinite topological endomorphism monoid. Moreover, restriction to any unbounded collection of free ranks preserves the module-valued functor category (Corollary~\ref{cor:cofinal-rank-invariance}).

When $A$ is finite, local, and self-injective, the full subcategory on finite $A$-submodules is inclusion-cofinal and is equivalent to the category $\C_A$ of finite $A$-modules. The Grothendieck topology on $\C_A^{\op}$ is then rigid with $\D_A$ as its rigid subcategory. Applying the Comparison Lemma one again gets the above equivalence. Moreover, $\C_A$ is a generalized Reedy category. Under the extra assumption that the characteristic of the residue field $A/\mathfrak m$ is invertible in $k$, we show that a standard module $\Delta_M$ is projective precisely when $M$ is semisimple. We then compute the morphism spaces and composition of a Morita-equivalent category, including an example over the dual numbers.

\subsection{Relation to previous works}

Rogers \cite{Rogers} constructed sites for continuous monoid actions using continuous principal actions and open congruences. In \cite{Rogers1}, he also interpreted continuous actions of topological endomorphism monoids of models through classifying topoi. Our construction provides a transformation realization adapted to finite uniform supports, making the finite-restriction reconstruction of \cite{LiRel} available for arbitrary topological
monoids. The criteria developed in this paper allow this general construction to be replaced by smaller, concrete models.

Many examples in this paper have been studied in the literature; see Magill's survey \cite{Magill1}. For instance: Magill \cite{Magill2} studied the factorization structure of the monoid of continuous self-maps on $\mathbb{R}$; Xu--Escard\'o \cite{XE} developed a constructive sheaf-theoretic model for the uniform continuity of functions on Cantor space; Elliott-Jonu\v{s}as-Mesyan-Mitchell-Morayne-P\'eresse \cite{EJMMMP} showed the compact-open topology on the monoid of continuous self-maps of Cantor space is the unique Polish semigroup topology, and Megrelishvili-Shlossberg \cite{MS} mentioned the importance of this monoid in the exploration of non-archimedean topological monoids. In this paper we impose a different topology on these monoids. Our purpose is also different: continuous actions on discrete sets and modules, and representations of associated finite-restriction categories.

On the representation-theoretic side, many combinatorial categories arising in the reconstruction have already been studied. For instance, the functor category over the category of finite free modules is familiar to experts in generic representation theory. For finite fields it appears in Kuhn's work \cite{Kuhn}; it has also been studied by Djament-Touz\'e \cite{DjamentTouze} over more general rings, while Putman-Sam \cite{PS} study closely related categories of finite free modules with split-injective morphisms over finite rings. We do not aim to replace these theories or to reprove their results. Instead, we supply a different viewpoint with respect to which familiar functor categories arise as continuous representation categories of topological monoids.

\subsection{Organization}

The paper is organized as follows. In Section 2 we establish the reconstruction method and its basic reductions. Various monoidal realizations of $\FA$-modules are given in Section 3, and in Section 4 we study infinite free modules over finite rings, their reconstruction sites, and the associated Morita model.

\subsection{Acknowledgement on the use of AI}

The research program, principal ideas and strategy, and the main mathematical constructions of this work were conceived by the authors. During the preparation of the paper, ChatGPT (OpenAI) was used as an auxiliary tool for exploratory mathematical discussion, assistance in locating relevant literature, checking intermediate arguments, and improving the exposition. All mathematical statements, proofs and bibliographic references in the final manuscript were independently written or checked by the authors, who take full responsibility for the content of the paper.

\section{A reconstruction method}

Throughout this section, let $\M$ be a topological monoid, and $k$ a nonzero commutative ring with identity. Given a nonempty set $\Omega$, denote by $\mathscr{F} = \Omega^{\Omega}$ the full transformation monoid equipped with the pointwise convergence topology (the product topology) with respect to the discrete topology on $\Omega$. For a submonoid $\E \leqslant \mathscr{F}$, we always equip $\E$ with the subspace topology. All actions are unital left actions and all representation spaces have the discrete topology.

\begin{definition}
A \emph{continuous $\M$-set} is a discrete set $X$ endowed with a continuous action $\M \times X \to X$. A \emph{continuous $k\M$-module} is a $k\M$-module whose underlying $\M$-set is continuous.
\end{definition}

Let $\Cont(\M)$ and $k\M \Mod^{\mathrm{cont}}$ be the corresponding categories. Since the target set is discrete, continuity of the action map is equivalent to continuity of every orbit map $m \mapsto m \cdot x$ for all $m \in \M$ and $x \in X$.

\begin{definition}
A \emph{left congruence} on $\M$ is an equivalence relation $R \subseteq \M^2$ such that
\[
(a,b) \in R \Longrightarrow (m \cdot a, m \cdot b) \in R
\]
for all $m \in \M$. It is \emph{open} if all its equivalence classes are open.
\end{definition}

Now we consider uniformly continuous representations of transformation monoids. In the case that $\E \leqslant \mathscr{F}$ is a transformation monoid equipped with the subspace topology, it is well known that $\E$ is a topological submonoid of $\mathscr{F}$.

\begin{definition}
Suppose that $\E \leqslant \mathscr{F}$. An $\E$-set $X$ is \emph{uniformly continuous} if for every $x \in X$ there is a finite subset $\Gamma \subseteq \Omega$ such that for any $f, g \in \E$, one has $f \cdot x = g \cdot x$ whenever $f|_{\Gamma} = g|_{\Gamma}$. In this case, we say that $\Gamma$ is a \emph{uniform support} of $x$.
\end{definition}

Denote the category of uniformly continuous $\E$-sets by $\UC(\E)$. One can similarly define uniformly continuous $k\E$-modules, and denote the category of these modules by $k\E \Mod^{\mathrm{uc}}$, which is an abelian category. It is clear from the definition that uniformly continuous $\E$-sets and $k\E$-modules are continuous; for details, see \cite{LiRel}.

Let $\Rcal_\M$ be the set of open left congruences on $\M$. For $R \in \Rcal_\M$, the quotient $\M/R$ is a continuous $\M$-set by left multiplication. Conversely, if $X$ is a continuous $\M$-set and $x \in X$, then
\[
R_x = \{(f, g) \in \M^2 \mid f \cdot x = g \cdot x \}
\]
is an open left congruence and the orbit $\M \cdot x$ is isomorphic to $\M/R_x$.

We now construct a transformation monoid $\E_\M$ over a discrete set $\Omega_{\M}$ such that continuous representations of $\M$ can be realized as uniformly continuous representations of $\E_\M$. Explicitly, let
\begin{equation}\label{eq:OmegaM}
\Omega_{\M} =\bigsqcup_{R \in \Rcal_\M} \M/R
\end{equation}
and equip $\Omega_\M$ with the discrete topology. Left multiplication induces a monoid homomorphism
\[
\eta: \M \to (\Omega_\M)^{\Omega_\M}.
\]
Let $\E_\M = \eta(\M)$ equipped with the subspace topology and the pointwise uniformity induced by $\Omega_M$.

\begin{theorem}\label{thm:universal-realization}
The restriction along $\eta: \M \twoheadrightarrow \E_\M$ induces equivalences
\[
\UC(\E_\M) \simeq \Cont(\M), \qquad k\E_\M \Mod^{\mathrm{uc}} \simeq k\M \Mod^{\mathrm{cont}}.
\]
Moreover, every continuous $\E_\M$-action is uniformly continuous.
\end{theorem}

\begin{proof}
We prove several facts, which together imply the desired conclusion.

\medskip
\noindent
\textbf{(1) The map $\eta$ is a continuous monoid homomorphism.} Given $R \in \Rcal_\M$, $\M$ acts on $\M/R$ by left multiplication: $a \cdot[b]_R=[ab]_R$, namely $\eta(a)([b]_R)=[ab]_R$. It follows that $\eta$ is a monoid homomorphism. It remains to prove continuity; that is, for every $\omega \in \Omega_\M$, the orbit map
\[
\M \longrightarrow \Omega_\M, \qquad a \longmapsto \eta(a)(\omega),
\]
is continuous. But this is clear. Indeed, write $\omega = [b]_R$ for some $R \in \Rcal_\M$. By the definition of $\Rcal_\M$, the action of $\M$ on $\M/R$ is continuous with respect to the discrete topology on $\M/R$, implying the continuity of the orbit map $a \mapsto a[b]_R=[ab]_R$.

\medskip
\noindent
\textbf{(2) A functor from $\Cont(\M)$ to $\UC(\E_\M)$.} Let $X$ be a continuous $\M$-set. For $x\in X$, since the action of $\M$ on $X$ is continuous, one has $R_x \in \Rcal_\M$, and so the quotient $\M/R_x$ is a component of $\Omega_M$.

We show that the $\M$-action on $X$ factors through $\E_\M$. Suppose that $\eta(a) = \eta(b)$. Evaluating these two transformations at the point $[1]_{R_x} \in \M/R_x$, we have
\[
[a]_{R_x} = \eta(a)([1]_{R_x}) = \eta(b)([1]_{R_x}) = [b]_{R_x},
\]
namely $a \cdot x = b \cdot x$. Since this holds for every $x \in X$, the action of $\M$ factors uniquely through the quotient monoid $\E_\M = \eta(\M)$. Thus $X$ is also an $\E_\M$-set. Moreover, the action of $\E_\M$ on $X$ is uniformly continuous since for each $x \in X$, the singleton $\Gamma_x = \{[1]_{R_x}\} \subseteq \Omega_\M$ is a uniform support for $x$. To see this, take $f, g \in \E_\M$ such that $f|_{\Gamma_x} = g|_{\Gamma_x}$, and choose $m, n \in \M$ such that $\eta(m) = f$ and $\eta(n) = g$. The equality on $\Gamma_x$ gives
\[
[m]_{R_x} = m([1]_{R_x}) = n([1]_{R_x}) = [n]_{R_x},
\]
namely $m \cdot x = n \cdot x$, and so $f \cdot x = g \cdot x$. Thus every $x \in X$ has a singleton uniform support, and hence the action of $\E_\M$ on $X$ is uniformly continuous.

We have shown that every continuous $\M$-set can be viewed as a uniformly continuous $\E_\M$-set. This construction is clearly functorial. In this way we obtain the desired functor.

\medskip
\noindent
\textbf{(3) The quasi-inverse functor.} Let $Y$ be a uniformly continuous $\E_\M$-set, which is automatically continuous. Since $\eta: \M \to \E_\M$ is continuous, restriction along $\eta$ gives a continuous $\M$-action on $Y$. Thus every uniformly continuous $\E_\M$-set can be viewed as a continuous $\M$-set. By this observation, we obtain a functor from $\UC(\E_\M)$ to $\Cont(\M)$.

The two functors are inverse to one another. Indeed, an $\M$-action factoring through $\eta: \M \twoheadrightarrow \E_\M$ has a unique induced $\E_\M$-action, while pulling that action back along $\eta$ recovers the original $\M$-action. Conversely, starting with an $\E_\M$-action and restricting along $\eta$, the factorization just constructed recovers the original $\E_\M$-action since $\eta$ is surjective. Moreover, the surjectivity of $\eta$ guarantees that morphisms are unchanged: given two uniformly continuous $\E_\M$-sets $X$ and $Y$, a map $f: X \to Y$ is $\M$-equivariant after restriction along $\eta$ if and only if $f(\eta(m) \cdot x) = \eta(m) \cdot f(x)$ for all $m \in \M$. Since every element of $\E_\M$ has the form $\eta(m)$, this is equivalent to $\E_\M$-equivariance.

In summary, we have shown that $\Cont(\M) \simeq \UC(\E_\M)$. The $k$-linearization version of this equivalence follows immediately.

\medskip
\noindent
\textbf{(4) Every continuous $\E_\M$-set is uniformly continuous.} Let $Y$ be a continuous $\E_\M$-set. As explained before, it can be viewed as a continuous $\M$-set via the restriction along $\eta$. For every $y \in Y$, $R_y$ is an open congruence, and hence $\M/R_y$ is a component of $\Omega_\M$. Exactly as above, the singleton $\{[1]_{R_y}\}$ is a uniform support for $y$, so $Y$ is a uniformly continuous $\E_\M$-set.
\end{proof}

\begin{remark}
Although the monoid homomorphism $\eta$ is continuous, the topology on $\M$ induced from the pointwise-convergence topology on $\E_\M$ may not coincide with the original topology on $\M$. In general, the induced topology is coarser. In other words, the inverse-image uniformity on $\M$ determined by $\eta$ records only the relations arising from continuous $\M$-actions. Its associated topology is the initial topology for all
orbit maps of continuous $\M$-sets, and hence is contained in the original topology of $\M$.

This distinction does not affect Theorem~\ref{thm:universal-realization}, where the equivalence is with respect to the \emph{original} topology on $\M$. Indeed, continuity of $\eta$ ensures that every uniformly continuous $\E_\M$-action restricts to a continuous $\M$-action, while the construction of $\Omega_\M$ ensures conversely that every continuous $\M$-action factors through $\E_\M$ and becomes uniformly continuous. Thus the theorem does not require $\eta$ to be a topological embedding, nor does it require the topology induced by $\eta$ to agree with the original topology on $\M$.
\end{remark}

The above theorem tells us that continuous representations of an arbitrary topological monoid $\M$ can be realized as uniformly continuous representations of a concrete transformation monoid $\E_\M$ over a discrete set $\Omega_{\M}$. In the rest of this section we briefly recall the finite-restriction reconstruction of \cite{LiRel} and the torsion-theoretic description of \cite{DLL}.

Define a small category $\C_\M$ as follows. Its objects are nonempty finite subsets $\Gamma \subseteq \Omega_\M$, and
\[
\C_\M(\Gamma, \Sigma) = \{f|_\Gamma: \Gamma \to \Sigma \mid f \in \E_{\M},\, f(\Gamma) \subseteq \Sigma\}.
\]
In particular, every inclusion of finite subsets is a morphism in $\C_\M$. We also define the following Grothendieck topology on $\C_\M^{\op}$: given a finite subset $\Gamma$ of $\Omega_\M$, a sieve $F \subseteq \C_\M^{\op}(-, \Gamma) = \C_\M(\Gamma, -)$ is a covering sieve if and only if it contains an inclusion $\iota: \Gamma \hookrightarrow \Sigma$ for some finite $\Sigma$.

The following theorem gives a monoidal version of Artin's reconstruction theorem.

\begin{theorem}[\cite{LiRel}] \label{thm:finite-reconstruction}
One has the following equivalence
\[
k\E_\M \Mod^{\mathrm{uc}} \simeq \Sh(\C_\M^{\op}, \, J, \, \underline{k}),
\]
where $\Sh(\C_\M^{\op}, \, J, \, \underline{k})$ is the category of sheaves of $k$-modules over the site $(\C_\M^{\op}, \, J)$.
\end{theorem}

An explicit pair $(\mathcal{L}, \mathcal{R})$ inducing the equivalence is described in \cite{LiRel}. For $V \in k\E_{\M} \Mod^{\mathrm{uc}}$ and a finite subset $\Gamma$ of $\Omega_\M$, define
\[
\mathcal L(V)(\Gamma) = \{v \in V \mid f \cdot v = g \cdot v \text{ whenever }f|_\Gamma = g|_\Gamma\}.
\]
If $u: \Gamma \to \Sigma$ is the restriction of $f \in \E_\M$, define $\mathcal L(V)(u)(v) = f \cdot v$. Uniform support makes this independent of the extension and the resulting functor $\mathcal{L}(V)$ is a sheaf of $k$-modules. Conversely, for a sheaf $F$, we take the filtered colimit over inclusions
\[
\mathcal R(F) = \varinjlim_{\Gamma \in \Fin(\Omega_\M)} F(\Gamma).
\]
If $[a]_\Gamma$ is represented by $a \in F(\Gamma)$, for $f \in \E_\M$, define
\[
f \cdot [a]_\Gamma = [F(f|_\Gamma)(a)]_{f(\Gamma)}.
\]
Then every element in $\mathcal R(F)$ has finite uniform support, and hence $\mathcal R(F)$ is a uniformly continuous $k\E_\M$-module.

Sheaves of $k$-modules over the site $(\C_\M^{\op}, J)$ have a torsion-theoretic interpretation, described in \cite{DLL}. Given a $\C_\M$-module $V$ (a presheaf of $k$-modules for $\C_\M^{\op}$) as well as an element $v \in V(\Gamma)$, we say that $v$ is \emph{$J$-torsion} if there is a covering sieve $S \subseteq \C_\M(\Gamma, -)$ such that $f \cdot v = 0$ for all $f \in S$. The $J$-torsion elements form a submodule $T_J(V)$. We say that $V$ is $J$-torsion if $T_J(V) = V$. Note that the full subcategory $\mathcal{T}_J$ consisting of $J$-torsion modules is a Serre localizing subcategory of $\C_\M \Mod$.

\begin{definition}[\cite{DLL}]
A $\C_\M$-module $V$ is \emph{$J$-saturated} if $T_J(V) = 0 = R^1T_J(V)$, where $R^1T_J$ is the first right derived functor of $T_J$. Equivalently, $V$ is $J$-saturated if
\[
\Hom_{k\C_\M}(T,V) = 0 = \Ext^1_{k\C_\M}(T,V)
\]
for every $T$ in $\mathcal{T}_J$.
\end{definition}

The following result gives a torsion-theoretic characterization of sheaves of $k$-modules.

\begin{theorem}[\cite{DLL}]\label{thm:torsion-sheaf}
Sheaves of $k$-modules over $(\C_\M^{\op}, J)$ are exactly $J$-saturated modules. Moreover, one has the following equivalence
\[
\Sh(\C^{\op}_{\M}, \,J, \, \underline{k}) \simeq \C_\M \Mod/\T_J,
\]
the Serre quotient of $\C_\M \Mod$ by the full subcategory of $J$-torsion modules.
\end{theorem}

Combining these equivalences together, we obtain the following result.

\begin{theorem}
Let $\M$ be a topological monoid. Then one has:
\[
k\M \Mod^{\mathrm{cont}} \simeq k\E_\M \Mod^{\mathrm{uc}} \simeq \Sh(\C_\M^{\op}, \, J, \, \underline{k}) \simeq \C_\M \Mod/\T_J.
\]
\end{theorem}

For a transformation monoid $\E \leqslant \Omega^\Omega$, the construction of $\Omega_\E$ is unnecessary if every continuous $\E$-action is uniformly continuous with respect to $\Omega$. The following proposition gives a sufficient condition.

\begin{definition}
Let $\Omega$ be a set. A family $\mathcal A$ of finite subsets of $\Omega$ is called \emph{inclusion-cofinal} if for every finite subset $\Gamma \subseteq \Omega$, there exists $\Sigma \in \mathcal A$ such that $\Gamma \subseteq \Sigma$.
\end{definition}

Given a transformation monoid $\E \leqslant \mathscr{F}$, an inclusion-cofinal family of finite subsets $\mathcal{A}$, as well as a certain $\Sigma \in \mathcal{A}$, we define a submonoid
\[
\E_\Sigma = \{ f \in \E \mid f|_\Sigma = \id_\Sigma \}.
\]
Those $\E_\Sigma$ with $\Sigma \in \mathcal{A}$ form a neighborhood basis of $\E$ at the identity.

\begin{proposition}\label{prop:common-right-multiple}
Let $\E$ be a submonoid of $\mathscr{F}$ equipped with the subspace topology. Suppose that for every finite $\Sigma$ in an inclusion-cofinal family $\mathcal{A}$ and every $f, g \in \E$ agreeing on $\Sigma$, there exist $h_1, h_2 \in \E_\Sigma$ such that $f \circ h_1 = g \circ h_2$. Then every continuous $\E$-set is uniformly continuous.
\end{proposition}

\begin{proof}
We show that each continuous $\E$-set $X$ is uniformly continuous. For every $x \in X$, since the orbit map
\[
\E \longrightarrow X, \qquad h \longmapsto h \cdot x
\]
is continuous at the identity, there exists a finite subset $\Gamma \subseteq \Omega$ such that $h \cdot x = x$ for every $h \in \E_\Gamma$. Choose a member $\Sigma \in \mathcal{A}$ such that $\Gamma \subseteq \Sigma$. If $f, g \in \E$ agree on $\Sigma$, then take $h_1, h_2 \in \E_\Sigma$ such that $f \circ h_1 = g \circ h_2$. Since $\E_\Sigma \subseteq \E_\Gamma$, it follows that both $h_1$ and $h_2$ fix $x$. Consequently, one gets
\[
f \cdot x = f \cdot(h_1 \cdot x) = (f \circ h_1) \cdot x = (g\circ h_2) \cdot x = g \cdot (h_2 \cdot x) = g \cdot x.
\]
Thus $\Sigma$ is a finite uniform support of $x$, and hence $X$ is a uniformly continuous $\E$-set.
\end{proof}

The second reduction has been described in \cite{LiRel}. We say that the category $\C_\M$ satisfies the \emph{finite retraction property} if for every inclusion morphism $\iota: \Gamma \hookrightarrow \Sigma$ in $\C_\M$, there is a morphism $p: \Sigma \to \Gamma$ such that $p \circ \iota = \id_\Gamma$. In that case, the Grothendieck topology $J$ is the trivial topology, and $\Sh(\C_\M^{\op}, \, J, \, \underline{k}) \simeq \C_\M \Mod$. The following result combines this observation with the Comparison Lemma.

\begin{proposition}\label{prop:cofinal-retractions}
Let $ \D \subseteq \C_\M$ be a full subcategory whose objects are inclusion-cofinal among finite subsets of $\Omega_\M$. Restriction induces an equivalence
\[
\Sh(\C_\M^{\op}, \, J, \, \underline{k}) \simeq \Sh(\D^{\op}, \, J_{\D}, \, \underline{k}).
\]
If every inclusion in $\D$ admits a retraction, then $J_{\D}$ is trivial and the right-hand category is $\D\Mod$.
\end{proposition}

\begin{proof}
Every object $\Gamma$ of $\C_\M$ is contained in an object $\Sigma$ of $\D$, so the opposite of the inclusion $\Gamma \hookrightarrow \Sigma$ is a covering morphism in $\C_\M^{\op}$. Thus $(\D^{\op}, J_{\D})$ is a full dense subsite, and the Comparison Lemma gives the first equivalence.

To establish the second statement, we observe that a covering inclusion $\Gamma \hookrightarrow \Sigma$ in $\C_\M$ with $\Gamma$ in $\D$ can be refined to an inclusion $\Gamma \hookrightarrow \Gamma'$ with $\Gamma'$ in $\D$. Consequently the induced topology $J_\D$ is generated by opposites of inclusions between objects of $\D$. Under the additional hypothesis these are split epimorphisms in $\D^{\op}$. Since every sieve containing a split epimorphism contains the identity and is maximal, the induced topology is trivial, and hence the sheaf category coincides with the presheaf category.
\end{proof}

The third reduction concerns the category $\C_\M$. For a few interesting transformation monoids, their associated categories $\C_\M$ have a generalized Reedy structure, so one can construct generalized Dold-Kan correspondences using the combinatorial machinery described in \cite{LiHom, LiNorm}.

We apply the reconstruction method and reductions described above to some well known examples.

\begin{example} \label{full example}
Let $\Omega$ be an infinite discrete set and $\mathscr{F} = \Omega^{\Omega}$. If $f, g \in \mathscr{F}$ agree on a nonempty finite subset $\Gamma$, one can choose a retraction $p: \Omega \to \Gamma$ such that $p|_\Gamma = \id_{\Gamma}$. Clearly, $p$ is contained in $\mathscr{F}_\Gamma$, and $f \circ p = g \circ p$. By Proposition \ref{prop:common-right-multiple}, every continuous $k\mathscr{F}$-module is uniformly continuous. It is easy to check that $\C_\mathscr{F}$ is the category of nonempty finite subsets of $\Omega$, which is equivalent to $\FA$. Furthermore, every inclusion in $\C_\mathscr{F}$ has a retraction. The above reductions then give \[
k\mathscr{F} \Mod^{\mathrm{cont}} = k\mathscr{F} \Mod^{\mathrm{uc}} \simeq \mathrm{FA} \Mod.
\]
\end{example}

\begin{example}
Let $\E = \End(\mathbb{Q}, \leqslant)$ be the monoid of order-preserving maps, where $\mathbb{Q}$ is given the discrete topology. If $f, g \in \E$ agree on a nonempty finite subset $\Gamma$, one can choose an order-preserving retraction $p: \mathbb{Q} \to \Gamma$ such that its restriction on $\Gamma$ is the identity map. Clearly, $p$ is contained in $\E_\Gamma$, and $f \circ p = g \circ p$. By Proposition \ref{prop:common-right-multiple}, every continuous $k\E$-module is uniformly continuous. By mimicking the argument in the previous example as well as Dold-Kan correspondence, we obtain
\[
k\E \Mod^{\mathrm{cont}} = k\E \Mod^{\mathrm{uc}} \simeq \Delta \Mod \simeq \mathrm{Ch}_{\geqslant 0} (k \Mod),
\]
where $\Delta$ is the simplex category.
\end{example}

\begin{example}
Let $\E = \End_{\mathbb{F}_q}(\mathbb{F}_q^{\infty})$, where $\mathbb{F}_q$ is a finite field, and
\[
\mathbb{F}_q^{\infty} = \varinjlim_n \mathbb{F}_q^n,
\]
equipped with the discrete topology. Again, applying Proposition \ref{prop:common-right-multiple} we can show that every continuous $k\E$-module is uniformly continuous. The other two reductions and Kuhn's theorem \cite{Kuhn} tell us
\[
k\E \Mod^{\mathrm{cont}} = k\E \Mod^{\mathrm{uc}} \simeq \mathrm{VA} \Mod \simeq \prod_{n \geqslant 0} k \mathrm{GL}_n(\mathbb{F}_q) \Mod
\]
when $q$ is invertible in $k$, where $\mathrm{VA}$ is the category of finite dimensional vector spaces over $\mathbb{F}_q$.
\end{example}

\begin{example}
Let $\E$ be the monoid of affine maps on $\mathbb{F}_q^{\infty}$. The same argument gives
\[
k\E \Mod^{\mathrm{cont}} = k\E \Mod^{\mathrm{uc}} \simeq \mathrm{AF} \Mod,
\]
where $\mathrm{AF}$ is the category of finite dimensional affine spaces over $\mathbb{F}_q$. Moreover, when $q$ is invertible in $k$, $\mathrm{AF} \Mod$ is equivalent to the representation category of a thin $k$-linear category. For details, see \cite{LiNorm}.
\end{example}

\section{Topological monoids realizing \texorpdfstring{$\FA$}{FA}--modules}

As we mentioned before, different topological monoids may have equivalent continuous representation categories. In this section we exhibit several such monoids whose continuous representation categories are equivalent to $\FA \Mod$. Throughout this section let $\Omega$ be an infinite set, and let $\mathscr F = \Omega^\Omega$ equipped with the product topology induced from the discrete topology on $\Omega$. All submonoids of $\mathscr{F}$ considered in this section are equipped with the corresponding subspace topology.

\begin{proposition}\label{prop:dense-full-FA}
Let $\E$ be a dense submonoid of $\mathscr F$. Then restriction induces an equivalence
\[
k\mathscr{E} \Mod^{\mathrm{uc}} \simeq k\mathscr{F} \Mod^{\mathrm{uc}} \simeq \FA \Mod.
\]
If every continuous $k\E$-module is uniformly continuous, then
\[
k\E\Mod^{\mathrm{cont}} \simeq \FA \Mod.
\]
\end{proposition}

\begin{proof}
The conclusion follows from the density criterion \cite[Corollary~7.12]{LiRel} and Example \ref{full example}.
\end{proof}

\begin{example}\label{ex:surjective-FA}
Let $\E = \operatorname{Sur}(\Omega)$ be the monoid of all surjective self-maps on $\Omega$. It is dense in $\mathscr F$ by \cite[Example~7.15]{LiRel}. We prove independently that every continuous $\E$-action is uniformly continuous. Given $f,g \in \E$ agreeing on a finite subset $\Gamma \subseteq \Omega$, form the pullback
\[
P = \{(a, b) \in \Omega \times \Omega \mid f(a) = g(b)\}.
\]
The coordinate projections $\pi_1,\pi_2:P\to\Omega$ are surjective. Since $\Omega$ is infinite, $|P| = |\Omega|$. We may therefore choose a surjection $r: \Omega \twoheadrightarrow P$ such that $r(\gamma) = (\gamma, \gamma)$ for all $\gamma \in \Gamma$. Then $h_i=\pi_i r$ are surjective self-maps fixing $\Gamma$ pointwise, and $f h_1=g h_2$. Thus Proposition~\ref{prop:common-right-multiple} applies and shows that every continuous $k\E$-module is uniformly continuous. Proposition~\ref{prop:dense-full-FA} now gives
\[
k\operatorname{Sur}(\Omega)\Mod^{\mathrm{cont}}   =k\operatorname{Sur}(\Omega)\Mod^{\mathrm{uc}} \simeq \FA \Mod.
\]
\end{example}

Let $\tau$ be a topology on $\Omega$, and let $\E \leqslant \mathscr{F}$ be a submonoid consisting of certain self-maps which are continuous with respect to $\tau$. We remind the reader that $\E$ always carries the subspace topology inherited from $\mathscr{F}$, which in general is different from the compact-open topology induced from the topology $\tau$ on $\Omega$.

We describe sufficient conditions such that the two assumptions in the previous proposition hold.

\begin{definition}\label{def:FEP}
The topological space $(\Omega, \tau)$ satisfies the \emph{finite extension property for self-maps} (\emph{FEP}) with respect to $\E$ if, for every finite subset $\Gamma \subseteq \Omega$ and every map $u: \Gamma \to \Omega$, there is a continuous self-map $\widetilde{u} \in \E$ such that
$\widetilde{u}|_{\Gamma} = u$.
\end{definition}

We describe an elementary condition equivalent to the first assumption in the previous proposition.

\begin{lemma}\label{lem:FEP-dense}
The submonoid $\E$ is dense in $\mathscr{F}$ if and only if $(\Omega, \tau)$ satisfies FEP with respect to $\E$.
\end{lemma}

\begin{proof}
We observe that the sets $ U(f, \Gamma) = \{g \in \mathscr{F} \mid g|_\Gamma = f|_\Gamma\}$ with $\Gamma$ a finite subset of $\Omega$ form a basis for the pointwise-convergence topology on $\mathscr{F}$. Consequently, $\E$ is dense in $\mathscr{F}$ if and only if every such neighborhood meets $\E$, and if and only if the restriction $f|_\Gamma$ extends to an element in $\E$.
\end{proof}

\begin{remark}
The reader can check that $\Omega$ satisfies FEP with respect to $\E$ if and only if the category of nonempty finite subsets of $\Omega$ and restrictions of elements in $\E$ coincides with the category of nonempty finite subsets of $\Omega$ and all maps.
\end{remark}

Here is a useful sufficient condition for the second hypothesis.

\begin{lemma}\label{lem:finite-retractions-FA}
Suppose that for every nonempty finite subset $\Gamma\subseteq\Omega$, there exists a continuous retraction $r_\Gamma\in\E$ from $\Omega$ onto $\Gamma$. Then every continuous $k\E$-module is uniformly continuous.
\end{lemma}

\begin{proof}
Given a continuous $k\E$-module $V$, we need to find a uniform support for any $v \in V$. Note that continuity at the identity of $\E$ gives a finite set $\Gamma$ such that every $h \in \E$ fixing $\Gamma$ pointwise also fixes $v$. If $\Gamma$ is empty, then all elements of $\E$ fix $v$; in particular, any nonempty finite subset $\Gamma$ of $\Omega$ serves as a uniform support of $v$. Otherwise, choose $r_\Gamma$ as in the assumption. For $f, g \in \E$ agreeing on $\Gamma$, we have $f \circ r_\Gamma = g\circ r_\Gamma$, $r_\Gamma v = v$, and hence
\[
  f \cdot v = (f \circ r_\Gamma) \cdot v = (g \circ r_\Gamma) \cdot v = g \cdot v.
\]
Thus $\Gamma$ is a uniform support of $v$. Consequently, $V$ is uniformly continuous.
\end{proof}

We apply the above results to a familiar example.

\begin{example}\label{ex:cantor-FA}
Let $\Omega$ be an infinite zero-dimensional Hausdorff space (for example, the classical Cantor space, and infinite Stone spaces), and $\E$ the monoid of all continuous self-maps on $\Omega$. For every finite set $\Gamma = \{x_1, \ldots, x_r\}$ of $\Omega$, we choose pairwise disjoint clopen subsets $U_1, \ldots, U_r$ with $x_i \in U_i$. Let $U_0 = \Omega \setminus \bigcup_{i=1}^r U_i$ and define
\[
  r_\Gamma(x)=
  \begin{cases}
    x_i, & x \in U_i, \quad 1 \leqslant i \leqslant r,\\
    x_1, & x \in U_0.
  \end{cases}
\]
This is a continuous retraction onto $\Gamma$, so every continuous $k\E$-module is uniformly continuous by Lemma \ref{lem:finite-retractions-FA}. Moreover, any map $u: \Gamma \to \Omega$ is continuous since the subspace topology on the finite set $\Gamma$ is the discrete topology. It follows that the composite $u \circ r_{\Gamma} \in \E$ is a continuous extension of $u$, so $\E$ is dense in $\mathscr{F}$ by Lemma \ref{lem:FEP-dense}. Consequently, Proposition \ref{prop:dense-full-FA} gives
\[
  k\E \Mod^{\mathrm{cont}} = k\E \Mod^{\mathrm{uc}} \simeq\FA\Mod.
\]
\end{example}

Let $\tau$ be the usual topology on the real line $\mathbb{R}$, and $\E$ be the monoid of continuous self-maps on $\mathbb{R}$, equipped with the subspace topology induced from $\mathscr{F}$. The real line satisfies FEP since any map from a finite subset of $\mathbb R$ to $\mathbb R$ extends by piecewise linear interpolation. Thus $\E$ is dense in $\mathscr{F}$. However, the preceding retraction argument does not apply to this example since a continuous image of $\mathbb R$ cannot be a finite set with more than one element. To show the uniform continuity of continuous $k\E$-modules, we introduce a new criterion.

\begin{lemma}\label{lem:common-image-interpolation}
Let $\Omega$ be an infinite topological space satisfying FEP, and $\E$ the monoid of all continuous self-maps on $\Omega$. Suppose that for every nonempty finite subset $\Gamma \subseteq \Omega$ and every $f, g \in \E$ with
$f|_\Gamma = g|_\Gamma$, there exists $p\in\E$ such that
\[
  p|_\Gamma=f|_\Gamma=g|_\Gamma,  \qquad p(\Omega) \subseteq f(\Omega) \cap g(\Omega).
\]
Then every continuous $k\E$-module is uniformly continuous.
\end{lemma}

\begin{proof}
Given a continuous $k\E$-module $V$ and $v \in V$, we need to find a uniform support for $v$. By the continuity of the orbit map at the identity, there is a finite subset $\Gamma \subseteq \Omega$ such that $h \cdot v = v$ for every $h \in \E_{\Gamma}$. If $\Gamma$ is empty, then $\E_{\Gamma} = \E$, namely every element of $\E$ fixes $v$. Consequently, any nonempty finite subset $\Gamma$ of $\Omega$ is a uniform support of $v$.

Now suppose that $\Gamma$ is nonempty, and $f, g \in \E$ agree on $\Gamma$. Choose $p$ as in the assumption. By continuity of the orbit map at $p$, there is a finite subset $\Sigma \subseteq \Omega$ such that
\begin{equation}\label{eq:common-image-local}
  q|_\Sigma = p|_\Sigma \quad \Longrightarrow \quad q \cdot v = p \cdot v, \qquad \forall \, q \in \E.
\end{equation}
Enlarging $\Sigma$ if necessary, we may assume that $\Gamma \subseteq \Sigma$. For each $x \in \Sigma \setminus \Gamma$, choose $b_x \in \Omega$ such that $f(b_x) = p(x)$, which is possible since $p(\Omega) \subseteq f(\Omega)$. By FEP, one can find a certain $h_f \in \E$ extending the map $\Sigma \to \Omega$ defined as follows:
\[
x \longmapsto
\begin{cases}
x, & x \in \Gamma,\\
b_x, & x \in \Sigma \setminus \Gamma.
\end{cases}
\]
Then $h_f \in \E_{\Gamma}$, and hence $h_f \cdot v = v$. Moreover, we get $(f \circ h_f)|_\Sigma = p|_\Sigma$: on $\Sigma \cap \Gamma$ this follows from $f|_\Gamma = p|_\Gamma$, and on $\Sigma \setminus \Gamma$ it follows from definition of $b_x$. By \eqref{eq:common-image-local} we obtain
\[
  f \cdot v = f \cdot (h_f \cdot v) = (f \circ h_f) \cdot v = p \cdot v.
\]
Similarly, one can carry out the same construction for $g$ to show $g \cdot v = p \cdot v$. Thus $f \cdot v = g \cdot v$, namely $\Gamma$ is actually a uniform support of $v$. Thus $V$ is uniformly continuous.
\end{proof}

\begin{proposition}\label{prop:ordered-space-FA}
Let $\Omega$ be an infinite linearly ordered set equipped with the order topology, and let $\E$ be the topological submonoid of $\mathscr{F}$ consisting of all continuous self-maps of $\Omega$. Suppose that $\Omega$ is connected and satisfies FEP. Then every continuous $k\E$-module is uniformly continuous, and
\[
k\E\Mod^{\mathrm{cont}} = k\E\Mod^{\mathrm{uc}} \simeq \FA \Mod.
\]
In particular, the conclusion holds if the group $\Homeo(\Omega)$ acts $2$-transitively on $\Omega$.
\end{proposition}

\begin{proof}
We first verify the hypothesis of Lemma \ref{lem:common-image-interpolation}. Let $f$ and $g$ be two elements in $\E$ such that they agree on a nonempty finite subset $\Gamma \subseteq \Omega$. Let $a = \min f(\Gamma) = \min g(\Gamma)$ and $b = \max f(\Gamma) = \max g(\Gamma)$. Since $\Omega$ is connected, the images $f(\Omega)$ and $g(\Omega)$ are convex, so $[a,b] \subseteq f(\Omega) \cap g(\Omega)$.

Define the clamping map
\[
c_{a,b}: \Omega \longrightarrow \Omega, \qquad
c_{a,b}(y)=
\begin{cases}
a, & y \leqslant a,\\
y, & a \leqslant y \leqslant b,\\
b, & y \geqslant b.
\end{cases}
\]
This is a continuous map. Consequently, $p = c_{a,b} \circ f$ belongs to $\E$ and satisfies
\[
p|_\Gamma=f|_\Gamma=g|_\Gamma, \qquad p(\Omega) \subseteq[a,b] \subseteq f(\Omega) \cap g(\Omega).
\]
By Lemma \ref{lem:common-image-interpolation}, every continuous $k\E$-module is uniformly continuous. Moreover, FEP implies that $\E$ is dense in $\mathscr{F}$. The asserted equivalence then follows from Proposition \ref{prop:dense-full-FA}.

It remains to show that $2$-transitivity of $\Homeo(\Omega)$ implies FEP. Let $\Gamma = \{x_1 < \cdots <x_n\}$ be a nonempty finite subset of $\Omega$, and let $u: \Gamma \to \Omega$ be a map. Write $y_i=u(x_i)$. If $n=1$, the constant map with value $y_1$ extends $u$, so we assume $n \geqslant 2$. For each $i$ with $y_i \neq y_{i+1}$, by the 2-transitivity, we can find $h_i \in \Homeo(\Omega)$ such that $h_i(x_i) = y_i$ and $h_i(x_{i+1}) = y_{i+1}$. Accordingly, we define a map $\widetilde{f}: \Omega \to \Omega$ by
\[
\widetilde{f}(x) =
\begin{cases}
y_1, & x \leqslant x_1,\\
h_i(x), & x_i \leqslant x \leqslant x_{i+1}  \text{ and } y_i \neq y_{i+1},\\
y_i, & x_i \leqslant x \leqslant x_{i+1}  \text{ and } y_i = y_{i+1},\\
y_n, & x \geqslant x_n.
\end{cases}
\]
It is easy to see that $\widetilde{f}$ is a continuous extension of $u$. Thus FEP holds.
\end{proof}

\begin{example}\label{ex:real-line-FA}
Let $\Omega = \mathbb{R}$ with its usual topology, and let $\E$ be the topological submonoid of $\mathscr{F}$ consisting of continuous self-maps of $\mathbb{R}$. The real line is connected, and its
homeomorphism group acts $2$-transitively. Thus by Proposition \ref{prop:ordered-space-FA}, we have
\[
k\E \Mod^{\mathrm{cont}} = k\E \Mod^{\mathrm{uc}} \simeq \FA \Mod.
\]

This example does not satisfy the common right multiple condition in Proposition \ref{prop:common-right-multiple}. Indeed, take a finite subset $\Gamma = \{a_1 < a_2 < \ldots < a_n \} \subseteq \mathbb{R}$ with $n \geqslant 2$, and define
\[
f(t) = \prod_{i = 1}^n (t - a_i)^2, \quad g(t) = -f(t).
\]
Clearly, $f|_\Gamma = g|_\Gamma$. Suppose there were continuous self-maps $h_1, h_2$ fixing $\Gamma$ pointwise and satisfying $f \circ h_1 = g \circ h_2$. Then $f \circ h_1$ must be the zero map, so $h_1$ must send $\mathbb{R}$ surjectively onto $\Gamma$, which is impossible.

The topology on $\E$ is important here. If instead we use the compact-open topology associated with the usual topology on $\mathbb R$, then $\E$ is path-connected. This forces every continuous orbit map from $\E$ to a discrete set to be constant, and hence every continuous action on a discrete set is trivial. Consequently, every continuous $k\E$-module has the trivial $\E$-action, and the continuous representation category is equivalent to $k \Mod$.
\end{example}

\section{\texorpdfstring{Endomorphism}{Endomorphism} monoids of infinite free modules over finite rings}\label{sec:direct-finite-ring}

In this section we describe another application of the method described in Section 2. Throughout this section let $A$ be an arbitrary finite ring, and
\[
\Omega = A^{\infty}=\varinjlim_r A^r
\]
with respect to the standard coordinate inclusions. Give $\Omega$ the discrete topology, and define $\E = \End_A(\Omega)$ with the pointwise--convergence topology.

\begin{lemma}\label{prop:auto-finite-ring}
Every continuous $k\E$-module is uniformly continuous.
\end{lemma}

\begin{proof}
The conclusion follows readily from Proposition~\ref{prop:common-right-multiple}. Indeed, the family $\mathcal{A} = \{ A^r \mid r \geqslant 0\}$ of standard coordinate submodules of $\Omega$ is inclusion-cofinal. Moreover, for $f, g \in \E$ agreeing on $A^r$, the standard coordinate projection $p_r: \Omega \to A^r$ satisfies $f \circ p_r = g \circ p_r$. Thus every continuous $k\E$-module is uniformly continuous.
\end{proof}

\subsection{The finite free site}

Denote by $\D$ the category whose objects are the finite free modules $A^r$, $r \geqslant 0$, and whose morphisms are all $A$-linear maps.

\begin{proposition}\label{finite-ring-reconstruction}
Let $k$ be a commutative ring. Then one has a natural equivalence
\[
k\E \Mod^{\mathrm{cont}} \simeq \D \Mod.
\]
\end{proposition}

\begin{proof}
By Lemma \ref{prop:auto-finite-ring}, every continuous representation of $\E$ is uniformly continuous. Since $A$ is finite, each coordinate submodule $A^n$ is a finite subset of $\Omega$. Moreover, every $A$-linear map $u: A^m \to A^n$ in $\D$ extends to an element $\iota_n u p_m$ in $\E$, where $p_m: \Omega \to A^m$ is the coordinate projection and $\iota_n: A^n \hookrightarrow \Omega$ is the inclusion. Thus $\D$ is a full subcategory of the finite-restriction category. Every inclusion in $\D$ splits by coordinate projections. The conclusion then follows from Theorem \ref{thm:finite-reconstruction} and Proposition \ref{prop:cofinal-retractions}.
\end{proof}

\begin{remark}
The conclusions of the previous proposition and Lemma \ref{prop:auto-finite-ring} actually hold for arbitrary unital rings. But here we do not pursue this full generality.
\end{remark}

For an unbounded subset $S \subseteq \mathbb N$, let $\D_S$ be the full subcategory of $\D$ consisting of objects $A^s$ with $s \in S$. The following result is very natural since every finite free module is a retract of a certain $A^s$ with $s \in S$.

\begin{corollary}\label{cor:cofinal-rank-invariance}
The restriction along the inclusion $i_S: \D_S \hookrightarrow \D$ induces an equivalence $i_S^*: \D \Mod \xrightarrow{\ \sim\ } \D_S \Mod$. In particular, every infinite full subcategory of $\D$ has an equivalent module category.
\end{corollary}

\begin{proof}
Since $S$ is unbounded, the family $\{A^s \mid s \in S \}$ remains cofinal among finite subsets of $\Omega$. Again, every $A$-linear map between two members in this family extends to an endomorphism of $\Omega$, and every standard inclusion between them has a left inverse. The equivalence then follows from Theorem \ref{thm:finite-reconstruction} and Proposition \ref{prop:cofinal-retractions}. Moreover, since both reconstruction functors are defined by restriction to the chosen finite pieces, the equivalence is naturally isomorphic to $i_S^*$.
\end{proof}

\subsection{The finite site}

We now replace the coordinate free modules by all finite $A$-submodules of $\Omega$ and describe the resulting site. Let $\C$ be the category whose objects are finite $A$-submodules of $\Omega$ and whose morphisms are restrictions of elements in $\E$. It is clear that objects of $\C$ are inclusion-cofinal among finite subsets of $\Omega$. Thus one may use $\C$ as well as the Grothendieck topology $J$ generated by inclusions of finite $A$-submodules of $\Omega$ to reconstruct continuous $k\E$-modules, namely
\[
k\E \Mod^{\mathrm{cont}} \simeq \Sh(\C^{\op}, \, J, \, \underline{k}).
\]

The category $\D$ is a full subcategory of $\C$: between coordinate free modules every $A$-linear map occurs as a restriction of an element of $\E$. For a general finite ring, however, $\C$ need not be equivalent to the category of all finite $A$-modules and all linear maps. Moreover, the Grothendieck topology $J$ in the above equivalence is not transparent.

From now on we assume that $A$ is a finite local self-injective ring with the maximal ideal $\mathfrak m$ and residue field $\mathbb{F}_q = A/\mathfrak m$, which by Wedderburn's theorem is a finite field of cardinality $q$.

\begin{lemma} \label{prop:C-all-modules}
The category $\C$ is equivalent to the category of finitely generated $A$-modules.
\end{lemma}

\begin{proof}
Since $A$ is a finite left self-injective ring, it is quasi-Frobenius. Hence the injective envelope $E(V)$ of every finite left $A$-module $V$ is finitely generated and projective (see \cite[\S 15]{Lam}). As $A$ is local, $E(V)$ is finite free. Thus every finite $A$-module embeds into $A^n \subseteq \Omega$ for some $n$, and hence is isomorphic to an object of $\C$. Since $\C$ and the category of finite $A$-modules have the same objects up to isomorphism, it suffices to show that every $A$-linear map between two objects of $\C$ is a morphism in $\C$. Given such an $A$-linear map $
f: U \to V$, choose $n$ such that $V \subseteq A^n \subseteq \Omega$. Since $A$ is self-injective, the finite free module $A^n$ is injective, so the composite $U \xrightarrow{\,f\,} V \hookrightarrow A^n$ extends across the inclusion $U \hookrightarrow \Omega$ to an $A$-linear map $\widetilde{f}: \Omega \to A^n$. Composing with the standard inclusion $A^n \hookrightarrow \Omega$ we obtain an endomorphism in $\End_A(\Omega)$ whose restriction to $U$ is $f$. Thus $f$ is indeed a morphism in $\C$.
\end{proof}

For a finite $A$-module $V$, let $\eta_V: V \hookrightarrow E(V)$ be its injective envelope. The next lemma identifies the Grothendieck topology $J$ in the above equivalence.

\begin{lemma}\label{thm:rigid-topology}
Given a finite $A$-module $V$, a sieve $S \subseteq \C(V, -)$ of $J$ on $\C^{\op}$ is a covering sieve if and only if it contains $\eta_V$. Consequently, $J$ is the rigid subcategory topology determined by the full subcategory $\D$ of finite free modules.
\end{lemma}

\begin{proof}
If $S \subseteq \C(V, -)$ is a covering sieve, then it contains an inclusion $\iota: V \to U$. Since $E(V)$ is injective, there is an $A$-module homomorphism $h: U \to E(V)$ such that $h \circ \iota = \eta_V$. Since $\iota \in S$, it follows that $\eta_V$ is contained in $S$ as well. Conversely, $\eta_V$ is itself an inclusion (up to composing with an isomorphism), so the sieve $S$ containing it is a covering sieve. This proves the first statement.

To show that $J$ is a rigid topology, we claim that $J$-irreducible objects in $\C$ are precisely finite free $A$-modules. Indeed, if $V$ is a finite free $A$-module, then $\eta_V$ is an isomorphism, so every sieve $S_V$ containing $\eta_V$ is the maximal sieve. Thus $V$ is $J$-irreducible. On the other hand, if $V$ is not finite free, then $\eta_V$ is a proper inclusion. Let $S_V$ be the principal sieve generated by $\eta_V$. If it is the maximal sieve, then we could find an $A$-module homomorphism $f: E(V) \to V$ such that $f \circ \eta_V$ is the identity on $V$. But this is impossible, since otherwise $\eta_V$ splits, and hence $V$ as a direct summand of $E(V)$ must be finite free as well. Consequently, the covering sieve $S_V$ is not the maximal sieve, and hence $V$ is not $J$-irreducible.

To finish the proof, it suffices to show that for any finite $A$-module $V$, the sieve
\[
T_V = \bigcup_{W \text{ irreducible}} \C(W, -) \circ \C(V, W)
\]
is a covering sieve. But this is clear. Indeed, since $E(V)$ is $J$-irreducible, and $\eta_V \in \C(V, E(V))$, it follows that $T_V$ contains $\eta_V$, and hence is a covering sieve.
\end{proof}

Combining these two lemmas, we obtain another proof of Proposition \ref{finite-ring-reconstruction} as follows. Since $J$ is a rigid topology, objects in $\T_J$ are $\C$-modules $V$ such that its values on finite free $A$-modules vanish. Thus one can use \cite[C2.2, Theorem 2.2.3]{Jo} to deduce the last equivalence in the following formula
\[
k\E \Mod^{\mathrm{cont}} \simeq \Sh(\C^{\op},\, J, \, \underline{k}) \simeq \C \Mod/\T_J \simeq \D \Mod,
\]
where the equivalence between sheaves and $\D$-modules is given by restriction, with quasi-inverse given by right Kan extension.

\subsection{Generalized Reedy structure}

Throughout this subsection assume again that $A$ is a finite, self-injective local ring with maximal ideal $\mathfrak m$ and residue field $\mathbb{F}_q$. Given a finite $A$-module $M$, we denote its length by $\ell(M)$. From now on we replace $\C$ by a skeleton of the category of finite $A$-modules. This does not alter its representation category. The category $\C$ has two wide subcategories, namely $\C^+$ consisting of monomorphisms, and $\C^-$ consisting of epimorphisms.

\begin{lemma}\label{prop:Reedy}
The triple $(\C^+, \C^-, \ell)$ gives a generalized Reedy structure over $\C$ in the sense of \cite{BM}.
\end{lemma}

\begin{proof}
We briefly check the four axioms. Firstly, every map $f: M \to N$ admits a factorization $M \twoheadrightarrow \operatorname{im}(f) \hookrightarrow N$, and this factorization is unique up to a unique isomorphism of the intermediate object. A noninvertible morphism in $\C^+$ strictly increases composition length and a noninvertible morphism in $\C^-$ strictly decreases it. The intersection of $\C^+$ and $\C^-$ is the groupoid of isomorphisms. Finally, if $f: M \twoheadrightarrow N$ and $\sigma \in \Aut_A(N)$ satisfy $\sigma \circ f = f$, then by the surjectivity of $f$ we have $\sigma = \id_N$.
\end{proof}

For a finite $A$-module $M$, we define $P_M = k\C(M,-)$. It has a submodule $\mathfrak{I}_M$ spanned by all non-injective $A$-module homomorphisms starting from $M$. Define the associated standard module $\Delta_M = P_M / \mathfrak{I}_M$.

In the rest of this subsection we assume that $q$ is invertible in $k$. When $A = \mathbb{F}_q$, the following result can be recovered from Kuhn's theorem.

\begin{theorem}\label{thm:regular-iff-semisimple}
Let $M$ be a finite $A$-module. Then $\Delta_M$ is projective if and only if $M$ is semisimple.
\end{theorem}

\begin{proof}
Suppose that $M \cong \mathbb{F}_q^m$ is semisimple. Let $\mathrm{VA}$ be the category whose objects are $\mathbb{F}_q^a$ with $a \in \mathbb{N}$ and whose morphisms are $\mathbb{F}_q$-linear maps, which can be identified with the full subcategory of $\C$ consisting of semisimple $A$-modules. Note that $\mathrm{VA}$ is also a generalized Reedy category, so we can define the corresponding ideal $\mathfrak{J}_M$ and standard module $\Delta_M' = Q_M / \mathfrak{J}_M$, where $Q_M = k\mathrm{VA}(M, -)$. Moreover, by Kov\'{a}cs' theorem \cite{Kovacs} and Kuhn's result \cite{Kuhn}, there is an idempotent $\Phi_M \in k\End_{\mathbb{F}_q}(M) = k\End_A(M)$ such that $Q_M \Phi_M \cong \Delta_M'$.

Thus it suffices to show that $P_M \Phi_M \cong \Delta_M$. But this follows from \cite[Definition 2.11 and Theorem 2.12]{LiNorm}. Indeed, by Kov\'{a}cs' construction, one has $\id_M - \Phi_M \in kD_M$, where $D_M \subseteq \End_A(M)$ is the set of non-invertible endomorphisms over $M$. Furthermore, given a non-invertible morphism $f: M \to N$ in $\C^-$ with source $M$, namely a proper surjection in $\C$, since $M$ is semisimple, so is $N$. Thus $f$ is also a non-invertible morphism in $\mathrm{VA}^-$, and hence $f \Phi_M = 0$ since $\Phi_M$ is a normalization idempotent for $\mathrm{VA}$. Consequently, $\Phi_M$ is indeed a normalization idempotent for $\C^-$.

Conversely, suppose that $M$ is not semisimple. Then $\mathfrak m M \neq0$, so one can choose a simple submodule $0 \neq K \subseteq \mathfrak mM$. Let $s: M \twoheadrightarrow M/K$ be the natural projection, which is an unfactorizable morphism in $\C^-$. The normalization fiber over $s$ (see \cite[Subsection 2.2]{LiNorm}) is
\[
F_{s,s} = \{\id_M + v \mid v \in \End_A(M), \, v(M) \subseteq K\}.
\]
Since $\mathfrak{m} K = 0$, it follows that $v(\mathfrak{m} M) = \mathfrak{m} v(M) = 0$. But $K \subseteq \mathfrak{m}M$, so $v(K) = 0$ and hence $v^2 = 0$. Consequently, $\id_M + v$ has an inverse $\id_M-v$ in $\Aut_A(M)$. Since every element of $F_{s,s}$ is an automorphism, $\Delta_M$ cannot be projective by \cite[Corollary 2.8]{LiNorm}. This finishes the proof.
\end{proof}

\begin{remark}
The nonprojectivity direction holds over every nonzero coefficient ring. The hypothesis $q$ is invertible in $k$ is needed only for projectivity of the standard modules $\Delta_M$ whenever $M$ is semisimple.
\end{remark}

Now we define a new family of projective generators for $\C \Mod$ and use it to construct a $k$-linear category Morita equivalent to $\C$ such that its morphism spaces have explicit representation theoretic meaning of $A$. Explicitly, define
\[
Q_M=
\begin{cases}
\Delta_M, & M\text{ semisimple},\\
P_M, & M \text{ nonsemisimple}.
\end{cases}
\]

\begin{lemma} \label{lem:Q-generators}
The family $\{Q_M\}$ is a small family of projective generators of $\C \Mod$.
\end{lemma}

\begin{proof}
We have shown that every $Q_M$ is projective. It suffices to show that each $P_M$ can be written as a finite direct sum of members in the new family. If $M$ is nonsemisimple, then $P_M = Q_M$, and the conclusion holds trivially. If $M \cong \mathbb{F}_q^b$ is semisimple, we can filter $P_M$ by the dimension of the image. Note that every image is also semisimple, and the associated graded piece of image dimension $d$ is
\[
k[\Sur_{\mathbb{F}_q} (\mathbb{F}_q^b, \, \mathbb{F}_q^d)] \otimes_{k\mathrm{GL}_d(\mathbb{F}_q)} \Delta_{\mathbb{F}_q^d},
\]
a finite direct sum of copies of $\Delta_{\mathbb{F}_q^d}$. Since these standard modules are projective, the filtration splits. Thus $P_M$ is a finite direct sum of copies of standard modules as well.
\end{proof}

Define a $k$-linear category $\K$, which shares the same objects as $\C$, and whose morphism spaces are
\[
\K(M,N) = \Hom_{k\C}(Q_N, Q_M)
\]
for any finite $A$-modules $M$ and $N$. By the categorical version of Morita equivalence, the functor
\[
\C \Mod \longrightarrow \K \Mod, \qquad V \longmapsto \bigl(M \mapsto \Hom_{k\C}(Q_M, V) \bigr)
\]
is an equivalence.

For a finite $A$-module $M$, write $\Top M = M /\mathfrak{m} M$ and $\Soc M = \{v \in M \mid \mathfrak{m} v = 0 \}$. Both are naturally vector spaces over $\mathbb{F}_q \cong A/\mathfrak{m}$.

\begin{proposition}\label{prop:four-blocks}
Let $M, N$ be finite $A$-modules. Then:
\begin{enumerate}
\item If $M \cong \mathbb{F}_q^b$ and $N \cong \mathbb{F}_q^d$, then
\[
\K(M, N) \cong
\begin{cases}
k\mathrm{GL}_b(\mathbb{F}_q),& b = d,\\
0, & b \neq d.
\end{cases}
\]
\item If $N \cong \mathbb{F}_q^d$ and $M$ is not semisimple, then
\[
\K(M, N) \cong k[\Sur_{\mathbb{F}_q} (\Top M, \, \mathbb{F}_q^d)].
\]
\item If $N$ is not semisimple and $M \cong \mathbb{F}_q^b$, then
\[
\K(M, N) \cong k[\Inj_A(\mathbb{F}_q^b, \, N)] \cong k[\Inj_{\mathbb{F}_q} (\mathbb{F}_q^b, \, \Soc N)].
\]
\item If both $M$ and $N$ are not semisimple, then
\[
\K(M, N) \cong k[\Hom_A(M, N)].
\]
\end{enumerate}
\end{proposition}

\begin{proof}
We show the conclusion case by case.

\medskip
\noindent
(1). As in the proof of Theorem \ref{thm:regular-iff-semisimple}, there are normalization idempotents $\Phi_M$ and $\Phi_N$ such that
\[
Q_M = \Delta_M \cong P_M \Phi_M, \quad Q_N = \Delta_N \cong P_N \Phi_N.
\]
Consequently, one has
\[
\K(M, N) \cong \Hom_{k\C} (P_N \Phi_N, P_M \Phi_M) \cong \Phi_N P_M \Phi_M = \Phi_N (k\C(M, N)) \Phi_M = \Phi_N (k\mathrm{VA} (M, N)) \Phi_M.
\]
The conclusion then follows from Kuhn's theorem.

\medskip
\noindent
(2). In this case, one has
\[
\K(M, N) = \Hom_{k\C} (Q_N, Q_M) \cong \Hom_{k\C}(P_N\Phi_N, P_M) \cong \Phi_N k[\Hom_A(M, N)] \cong \Phi_N k[\Hom_{\mathbb F_q}(\Top M, N)].
\]
By Kuhn's decomposition theorem \cite{Kuhn}, the normalization idempotent $\Phi_N$ annihilates the contribution of all non-surjective maps, and its image is naturally identified with the free $k$-module on the surjective maps. Consequently,
\[
\Phi_N k[\Hom_{\mathbb F_q}(\Top M,\mathbb F_q^d)] \cong k[\Sur_{\mathbb F_q}(\Top M,\mathbb F_q^d)],
\]
and $\K(M,N) \cong k[\Sur_{\mathbb F_q}(\Top M,\mathbb F_q^d)]$ as claimed.

\medskip
\noindent
(3). This follows from Yoneda's lemma:
\[
\K(M, N) \cong \Hom_{k\C}(P_N, \Delta_{\mathbb{F}_q^b}) \cong \Delta_{\mathbb{F}_q^b}(N) \cong k[\Inj_A(\mathbb{F}_q^b,N)] \cong k[\Inj_{\mathbb{F}_q} (\mathbb{F}_q^b, \, \Soc N)],
\]
since every $A$-linear map from $\mathbb{F}_q^b$ to $N$ has image contained in $\Soc N$, and the $A$-module structure on $\mathbb{F}_q^b$ and $\Soc N$ both factor through $A/\mathfrak{m} \cong \mathbb{F}_q$.

\medskip
\noindent
(4). This follows from Yoneda's lemma.
\end{proof}

\begin{remark}
Set $e_M = \Phi_M$ for semisimple $M$, and $e_M = \id_M$ otherwise. Then
\[
\K(M,N) \cong e_N\,k[\Hom_A(M,N)]\,e_M.
\]
The identity at object $M$ in $\mathscr{K}$ is $e_M$, and composition is the bilinear extension of composition in $\C$:
\[
(e_L g e_N) \circ_{\K} (e_N f e_M) = e_L g e_N f e_M.
\]
Although the idempotents $\Phi_M$ may depend on choices, different choices yield equivalent projective generators $Q_M \cong \Delta_M$ and hence equivalent Morita models.
\end{remark}

\begin{example}\label{ex:dual-numbers}
Let $A=\mathbb{F}_q[\varepsilon]/(\varepsilon^2)$, $L = A/(\varepsilon)$, and assume $q$ is invertible in $k$. Then every finite $A$-module has the form $M_{a,b} = A^a \oplus L^b$, whose top and socle both have dimension $a+b$ over $\mathbb F_q$. Theorem \ref{thm:regular-iff-semisimple} asserts that $\Delta_{M_{a,b}}$ is projective exactly when $a = 0$.

Let us do a concrete composition for morphisms between $A$ and $L$ in $\mathscr{K}$. In this case, $e_L = \Phi_L = [\id_L]-[0_L]$, where $0_L$ is the zero endomorphism from $L$ to itself. For each $\lambda \in \mathbb{F}_q^{\times}$, define
\[
p_{\lambda}: A \longrightarrow L, \quad x \longmapsto \lambda \bar{x}
\]
and
\[
j_{\lambda}: L \longrightarrow A, \quad x \longmapsto \varepsilon \lambda x.
\]
By Proposition \ref{prop:four-blocks}, $\mathscr K(A,L)$ is a free $k$-module with basis
\[
\beta_\lambda = e_L[p_\lambda] = [p_\lambda] - [0_{A,L}], \qquad \lambda\in \mathbb{F}_q^\times.
\]
Similarly, $\mathscr K(L,A)$ is a free $k$-module with basis
\[
\alpha_\mu = [j_\mu] e_L = [j_\mu] - [0_{L,A}], \qquad \mu \in\mathbb{F}_q^\times.
\]
Since $p_\lambda j_\mu = 0_L$, whereas $j_\mu p_\lambda$ is multiplication by $\varepsilon\lambda\mu$, the composition rule gives
\[
\beta_\lambda \circ_{\mathscr{K}} \alpha_\mu = 0 \in \mathscr K(L,L),
\]
while
\[
\alpha_\mu \circ_{\mathscr{K}} \beta_\lambda = [m_{\varepsilon\lambda\mu}] - [0_A] \in \mathscr{K}(A,A).
\]
The latter element is nonzero because its two terms are distinct basis vectors.
\end{example}

\end{document}